\documentclass[12pt]{amsart}

\usepackage{amsmath,xspace,amssymb,mathrsfs}
\usepackage{color}

\input xy
\xyoption{all}
\xyoption{2cell}
\UseAllTwocells
\CompileMatrices

\renewcommand{\phi}{\varphi}

\newcommand{\Ker}{\operatorname{Ker}}
\newcommand{\Coker}{\operatorname{Coker}}
\newcommand{\Ima}{\operatorname{Im}}

\newcommand{\Ann}{\operatorname{Ann}}
\newcommand{\Hom}{\operatorname{Hom}}
\newcommand{\colimit}{\operatorname{colim}}
\newcommand{\Identity}{\operatorname{Id}}

\newtheorem{proposition}{Proposition}[section]
\newtheorem{lemma}[proposition]{Lemma}

\newtheorem{corollary}[proposition]{Corollary}
\newtheorem{theorem}[proposition]{Theorem}

\theoremstyle{definition}
\newtheorem{definition}[proposition]{Definition}
\newtheorem{example}[proposition]{Example}

\newtheorem{remark}[proposition]{Remark}

\begin{document}

\title{Gabriel localization theory and its applications}

\author{Abolfazl Tarizadeh}
\address{Department of Mathematics, Faculty of Basic Sciences, University of Maragheh \\
P. O. Box 55181-83111, Maragheh, Iran.
 }
\email{ebulfez1978@gmail.com}

\date{}

\begin{abstract} In this article first we develop the Gabriel localizations (abbreviated as G-localizations) for commutative rings, specially some new results in this direction are proven. Then, as an application, it is shown that a ring map is a flat epimorphism if and only if it corresponds to a kind of the G-localizations. As a by-product of this study, a characterization for the flatness of the quotient rings is given. The exactness of the G-localization functor are characterized. The structure of prime ideals in the G-localization rings are also studied. Finally, it is shown that the Gabriel localization theory is a natural generalization of the usual localization theory.
\end{abstract}

\maketitle

\section{Introduction}
Our main aim in this article is to analyze and then understand the structure of  flat epimorphisms of rings more deeply. Indeed, we have successfully applied the Gabriel localization theory and thereby we have found a new and simple proof to the fact that ``a ring map is a flat epimorphism if and only if it corresponds to a type of the G-localizations", see Corollary \ref{G-coro 2} and Theorem \ref{remark 200}. We should mention that, according to the knowledge of the author, there are some different proofs of this theorem or a part of it in the literature. E.g. see  \cite{Akiba}, \cite{Nastasescu-Popescu}, \cite{Olivier},  \cite{Popescu-Spircu} and \cite[Chap. XI, Theorem 2.1]{Stenstrom}. These proofs are based upon other enormous results and therefore it is practically very hard to follow and fully  understand the proof. In fact, they use the theory of Grothendieck categories and Giraud subcategories and also they use the fact that ``the set of idempotent topologizing systems on the ring $R$ is bijectively corresponding to the set of Giraud subcategories of the category of $R-$modules", see \cite[Chap. X, Theorem 2.1]{Stenstrom}. While the latter in turn is a huge result. Our approach, unlike their methods, is based upon some simple observations. More precisely, after developing the Gabriel localization theory for commutative rings we then obtain more general and new results,  namely Theorems \ref{G-lemma 10}, \ref{G-th 5} and \ref{G-lemma 1}. These results pave the way in order to reach to a natural and simple proof of the fact. Using these results then we are also able to characterize the exactness of the G-localization functor, see Theorem \ref{th621890}. Note that the G-localization functor, in general, is left exact. Moreover, as a by-product of this study, a characterization for the flatness of the quotient rings is given, Corollary \ref{coro 87}. \\

The localization theory, in particular the local rings, play a major role in commutative algebra, algebraic and arithmetic geometry, number and valuation theories. The notion of the G-localization with respect to an idempotent topologizing system, which is in turn a natural generalization of the usual localization theory, first appeared in the Gabriel thesis \cite[Chap V, \S 2]{Gabriel} which was conducted under the supervision of Grothendieck. In the literature, an idempotent topologizing system is also called a Gabriel topology. The Gabriel localization theory provides a very general method of localization which is even applicable in noncommutative situations. In commutative algebra a large number of important constructions are special cases of the G-localizations. The idea of the G-localization is essentially due to Grothendieck and it appeared in more general setting in SGA 4, tomme 1, expos\'{e} II. In this article, our presentation of the theory will closely follow the Gabriel's approach in Bourbaki \cite[Les exercices 16 \`{a} 25]{Bourbaki}.\\

The article is organized as follows. In Section 2, as mentioned in the above, the Gabriel localization theory is developed for commutative rings. The content of the third section was completely described in the first paragraph. Regarding to this section, we should mention that this expression of flat epimorphisms
(Corollary \ref{G-coro 2}) has some important applications. For example, using this result then one can show that every injective flat epimorphism of rings which is also of finite type then it is of finite presentation. The latter is also a highly non-trivial result in commutative algebra. In the final section, in addition to the characterizing the exactness of the G-localization functor,
the structure of prime ideals in the G-localization rings are also studied. Theorems \ref{th621890} and \ref{theorem 33} are the main results of this section.\\

\section{G-localizations}

Throughout the article, all of the rings which we consider are commutative.\\

\begin{definition} A \emph{topologizing system} on the ring $R$ is a non-empty family $\mathscr{F}$ of ideals of $R$ satisfying in the following conditions.\\
$\textbf{(a)}$ Every ideal of $R$ containing an ideal $I\in\mathscr{F}$ belongs to $\mathscr{F}$.\\
$\textbf{(b)}$ The family $\mathscr{F}$ is stable under finite intersections.\\
\end{definition}

\begin{example}\label{exam 1} Let $R$ be a ring, let $S$ be a multiplicative subset of $R$ and let $\mathcal{P}$ be a family of prime ideals of $R$. Examples of topologizing systems are the single-point set $\{R\}$, the set of ideals of $R$ containing a fixed ideal, the set of ideals $I$ of $R$ such that $I+J=R$ where $J$ is a fixed ideal of $R$, the set of ideals of $R$ which meeting $S$, the set of ideals $I$ of $R$ such that $\mathcal{P}\cap V(I)=\emptyset$, the set of ideals $I$ of $R$ such that $V(I)$ is contained in $V(J)$ where $J$ is a fixed ideal of $R$. See also Theorem \ref{G-lemma 1} and \cite[Chap. X, Theorem 2.1]{Stenstrom}. By $V(I)$ we mean the set of prime ideals $\mathfrak{p}$ of $R$ such that $I\subseteq\mathfrak{p}$. \\
\end{example}

\begin{definition} Let $\mathscr{F}$ be a topologizing system on the ring $R$. An $R-$module $M$ is called $\mathscr{F}-$\emph{negligible} if the annihilator of every element of $M$ belongs to $\mathscr{F}$. Clearly submodules, quotients, localizations and finite direct sums of $\mathscr{F}-$negligible modules are $\mathscr{F}-$negligible.\\
\end{definition}

Let $\mathscr{F}$ be a topologizing system on the ring $R$. For every $R-$module $M$, the set $\mathscr{F}(M)=\{m\in M : \Ann_{R}(m)\in\mathscr{F}\}$ is a $R-$submodule of $M$ because for every elements $m,m'\in\mathscr{F}(M)$ and for each $r\in R$ we have $\Ann(m)\cap\Ann(m')\subseteq\Ann(m+m')$ and $\Ann(m)\subseteq\Ann(rm)$. Clearly it is the greatest $\mathscr{F}-$negligible submodule of $M$.
Each $R-$linear map $u:M\rightarrow N$ induces a map $\mathscr{F}(M)\rightarrow\mathscr{F}(N)$ given by $m\rightsquigarrow u(m)$ which we denote it by $\mathscr{F}(u)$. In fact $\mathscr{F}(-)$ is a left exact functor from the category of $R-$modules to itself. \\

Let $\mathscr{F}$ be a topologizing system on the ring $R$. The family $\mathscr{F}$ with the relation $I<J$ if $J$ is a proper subset of $I$, is a directed poset. If $I\leq J$ then for each $R-$module $M$ the canonical injection $J\subseteq I$ induces the $R-$linear map $u_{I,J}:\Hom_{R}(I,M)\rightarrow\Hom(J,M)$ given by $f\rightsquigarrow f|_{J}$. Clearly $\big(\Hom_{R}(I,M), u_{I,J}\big)$ is an inductive system of $R-$modules and $R-$homomorphisms over the directed poset $(\mathscr{F},<)$. We shall denote by $M_{(\mathscr{F})}$ the inductive limit (colimit) of the system. Therefore $M_{(\mathscr{F})}=\colimit_{I\in\mathscr{F}}\Hom_{R}(I,M)$.
For each $I\in\mathscr{F}$ the canonical map $\Hom_{R}(I,M)\rightarrow M_{(\mathscr{F})}$, if there is no confusion, is denoted by $[\:]$. The composed map $\xymatrix{M\ar[r]^{\simeq\:\:\:\:\:\:\:\:\:\:\:\:\:\:\:}&
\Hom_{R}(R,M)\ar[r]&M_{(\mathscr{F})}}$ is also denoted by $\delta_{M}$. Clearly $\Ker(\delta_{M})=\mathscr{F}(M)$. Moreover $\Coker(\delta_{M})$ is $\mathscr{F}-$negligible. Because let $f:I\rightarrow M$ be a $R-$linear map   where $I\in\mathscr{F}$. For each $a\in I$, $(\delta_{m})|_{I}=a.f$ where $m=f(a)$ and $\delta_{m}: R\rightarrow M$ is given by $r\rightsquigarrow rm$. This means that $I\subseteq\Ann_{R}\big([f]+\Ima(\delta_{M})\big)$.\\

Every $R-$linear map $u:M\rightarrow N$, by the universal property of the colimits, induces a unique $R-$linear map $u_{(\mathscr{F})}:M_{(\mathscr{F})}\rightarrow N_{(\mathscr{F})}$ such that for each $I\in\mathscr{F}$ the following diagram is commutative $$\xymatrix{
\Hom_{R}(I,M)\ar[r]^{} \ar[d]^{} & \Hom_{R}(I,N)\ar[d]^{} \\ M_{(\mathscr{F})}\ar[r]^{\:\:\:\:u_{(\mathscr{F})}} & N_{(\mathscr{F})}}$$ where the columns are the canonical maps and the top row map is given by $f\rightsquigarrow u\circ f$. Therefore $u_{(\mathscr{F})}$ is defined as $[f]\rightsquigarrow[u\circ f]$. \\

\begin{lemma}\label{G-lemma 4} Let $\mathscr{F}$ be a topologizing system on the ring $R$. If \\ $\xymatrix{0\ar[r]&M'\ar[r]^{u}&M\ar[r]^{v}&M''}$ is an exact sequence of $R-$modules then the sequence $\xymatrix{0\ar[r]&M'_{(\mathscr{F})}\ar[r]^{u_{(\mathscr{F})}}&
M_{(\mathscr{F})}\ar[r]^{v_{(\mathscr{F})}}&M''_{(\mathscr{F})}}$ is exact.\\
\end{lemma}

{\bf Proof.} It is an easy exercise. $\Box$ \\

Let $\mathscr{F}$ and $\mathscr{G}$ be two topologizing systems on the ring $R$. We denote by $\mathscr{F}.\mathscr{G}$ the set of ideals $I$ of $R$ such that there exists some $J\in\mathscr{G}$ containing $I$ in which $J/I$ is $\mathscr{F}-$negligible. Clearly $\mathscr{F}.\mathscr{G}$ is a topologizing system on the ring $R$. Also $IJ\in\mathscr{F}.\mathscr{G}$ for all $I\in\mathscr{F}$ and all $J\in\mathscr{G}$. In particular, $\mathscr{F}$ and $\mathscr{G}$ are contained in $\mathscr{F}.\mathscr{G}$. A topologizing system $\mathscr{F}$ is called \emph{idempotent} if $\mathscr{F}.\mathscr{F}=\mathscr{F}$. All of the topologizing systems in Example \ref{exam 1}, except the second one, are idempotent.\\

\begin{proposition} Let $\mathscr{F}$ and $\mathscr{G}$ be two topologizing systems on the ring $R$ and let $M$ be a $R-$module. Then $M$ is $\mathscr{F}.\mathscr{G}-$negligible if and only if there exists a  $\mathscr{F}-$negligible submodule $M'$ of $M$ such that $M/M'$ is $\mathscr{G}-$negligible. In particular, if $\mathscr{H}$ is a third topologizing system on the ring $R$ then we have $(\mathscr{F}.\mathscr{G}).\mathscr{H}=
\mathscr{F}.(\mathscr{G}.\mathscr{H})$.\\
\end{proposition}

{\bf Proof.} If $M$ is $\mathscr{F}.\mathscr{G}-$negligible then take $M'=\mathscr{F}(M)$ and the remaining assertions are straightforward. The converse is also routine. $\Box$ \\

\begin{lemma}\label{G-lemma 5} Let $R$ be a ring and let $\mathscr{F}$ be a non-empty family of ideals of $R$. Then $\mathscr{F}$ is an idempotent topologizing system on the ring $R$ if and only if it satisfies in the following conditions.\\
$\textbf{(i)}$ Every ideal of $R$ containing an ideal $I\in\mathscr{F}$ belongs to $\mathscr{F}$.\\
$\textbf{(ii)}$ If $I\in\mathscr{F}$ and $J$ is an ideal of $R$ such that $J:a\in\mathscr{F}$ for all $a\in I$, then $J\in\mathscr{F}$.\\
\end{lemma}

{\bf Proof.} Suppose $\mathscr{F}$ is an idempotent topologizing system on the ring $R$. Let $I\in\mathscr{F}$ and let $J$ be an ideal of $R$ such that $J:a\in\mathscr{F}$ for all $a\in I$. It suffices to show that $J\in\mathscr{F}.\mathscr{F}$. But it is clear since $J'=I+J\in\mathscr{F}$ and for each $a\in I$, $\Ann_{R}(a+J)=J:a\in\mathscr{F}$. \\Conversely, suppose $I,J\in\mathscr{F}$. We show that $I\cap J\in\mathscr{F}$. But $J\subseteq IJ:a$ for all $a\in I$. Therefore $IJ\in\mathscr{F}$ and so $I\cap J\in\mathscr{F}$. Hence, $\mathscr{F}$ is a topologizing system. Take $J\in\mathscr{F}.\mathscr{F}$. Thus there exists some $I\in\mathscr{F}$ containing $J$ such that $I/J$ is $\mathscr{F}-$negligible. Therefore $J:a\in\mathscr{F}$ for all $a\in I$. Thus $J\in\mathscr{F}$. $\Box$ \\

From now onwards, if it is not stated, $\mathscr{F}$ is always an idempotent topologizing system on the ring $R$.\\

Here some new rings and modules are introduced. The basic set-up is as follows. Let $M$ be a $R-$module.
For each $I,J\in\mathscr{F}$ and for each $f\in\Hom_{R}(I,R)$ we have $f^{-1}(J)\in\mathscr{F}$. Because for each $a\in I$, $J\subseteq\Ann_{R}\big(a+ f^{-1}(J)\big)$. Therefore $I/f^{-1}(J)$ is $\mathscr{F}-$negligible and so
$f^{-1}(J)\in\mathscr{F}.\mathscr{F}=\mathscr{F}$. The map $f$ induces a map $f^{-1}(J)\rightarrow J$ which we denote it as usual by $f|_{f^{-1}(J)}$. Now we define the pairing $R_{(\mathscr{F})}\times M_{(\mathscr{F})}\rightarrow M_{(\mathscr{F})}$ as
$[f].[g]=[g\circ\big(f|_{f^{-1}(J)}\big)]$ where $J$ is the domain of $g$. Note that it is well-defined. Because let $[f_{1}]=[f_{2}]$ and $[g_{1}]=[g_{2}]$ where $f_{k}\in\Hom_{R}(I_{k},R)$ and $g_{k}\in\Hom_{R}(J_{k},M)$ for some $I_{k},J_{k}\in\mathscr{F}$ with $k=1,2$. There exist $L',L''\in\mathscr{F}$ such that $I_{k}\leq L'$, $J_{k}\leq L''$, $(f_{1})|_{L'}=(f_{2})|_{L'}$ and $(g_{1})|_{L''}=(g_{2})|_{L''}$. Take $L=f_{1}^{-1}(L'')\cap f_{2}^{-1}(L'')\cap L'$ which belongs to $\mathscr{F}$ and we have $\big(g_{1}\circ (f_{1})|_{f_{1}^{-1}(J_{1})}\big)|_{L}=\big(g_{2}\circ (f_{2})|_{f_{_{2}}^{-1}(J_{2})}\big)|_{L}$. The pairing is also $R-$bilinear (details omitted). \\

 In particular, the binary operation $R_{(\mathscr{F})}\times R_{(\mathscr{F})}\rightarrow R_{(\mathscr{F})}$, as multiplicative, puts a commutative ring structure on the $R-$module $R_{(\mathscr{F})}$. Its commutativity implies from this simple observation that for every $R-$linear maps $f,g:I\rightarrow R$ where $I$ is an ideal of $R$ and for every elements $a,a'\in I$ then $g\big(f(aa')\big)=f\big(g(aa')\big)$. The canonical map $\delta_{R}:R\rightarrow
R_{(\mathscr{F})}$ is a ring homomorphism. Moreover, the map
$R_{(\mathscr{F})}\times M_{(\mathscr{F})}\rightarrow M_{(\mathscr{F})}$ puts a $R_{(\mathscr{F})}-$module structure over $M_{(\mathscr{F})}$. For every $R-$linear map $u:M\rightarrow N$ then $u_{(\mathscr{F})}:M_{(\mathscr{F})}\rightarrow N_{(\mathscr{F})}$ is $R_{(\mathscr{F})}-$linear. Indeed, $M\rightsquigarrow M_{(\mathscr{F})}$ is a left exact functor from the category of $R-$modules to the category of $R_{(\mathscr{F})}-$modules.\\

\begin{definition} Let $M$ be a $R-$module. The $R_{(\mathscr{F})}-$module $$\big(M/\mathscr{F}(M)\big)_{(\mathscr{F})}$$ is called the \emph{Gabriel localization} of $M$ with respect to the system $\mathscr{F}$ and it is denoted by $M_{\mathscr{F}}$. Therefore $M_{\mathscr{F}}=
\colimit_{I\in\mathscr{F}}\Hom_{R}\big(I,M/\mathscr{F}(M)\big)$.
The composed map $\xymatrix{M\ar[r]^{\pi\:\:\:\:\:\:\:\:\:\:\:}&
M/\mathscr{F}(M)\ar[r]^{\:\:\:\:\:\:\:\delta}
&M_{\mathscr{F}}}$ is denoted by $j_{M}$. Therefore for each $m\in M$, $j_{M}(m)=[\delta_{\overline{m}}]$ where $\overline{m}=m+\mathscr{F}(M)$ and $\delta_{\overline{m}}:R\rightarrow M/\mathscr{F}(M)$ is defined as $r\rightsquigarrow rm+\mathscr{F}(M)$.\\
\end{definition}

\begin{lemma}\label{G-lemma 3} Let $M$ be a $R-$module. Then $\Ker(j_{M})=\mathscr{F}(M)$ and $\Coker(j_{M})$ is $\mathscr{F}-$negligible.\\
\end{lemma}

{\bf Proof.} First note that $\mathscr{F}\big(M/\mathscr{F}(M)\big)=0$ because if $\overline{m}=m+\mathscr{F}(M)\in\mathscr{F}\big(M/\mathscr{F}(M)\big)$ then $\Ann(m):a=\Ann(am)\in\mathscr{F}$ for all $a\in\Ann_{R}(\overline{m})$. Therefore, by Lemma \ref{G-lemma 5}, $\Ann(m)\in\mathscr{F}$ and so $\overline{m}=0$. We have $\Ker(j_{M})=j_{M}^{-1}(0)=\pi^{-1}\big(\Ker(\delta)\big)=
\pi^{-1}\Big(\mathscr{F}\big(M/\mathscr{F}(M)\big)\Big)=
\Ker(\pi)=\mathscr{F}(M)$. We also have $\Coker(j_{M})=\Coker(\delta)$ therefore $\Coker(j_{M})$ is $\mathscr{F}-$negligible. $\Box$ \\

\begin{lemma}\label{G-lemma 7} Let $M$ be a $R-$module. Then $M$ is $\mathscr{F}-$negligible if and only if  $M_{\mathscr{F}}=0$.\\
\end{lemma}

{\bf Proof.} If $M$ is $\mathscr{F}-$negligible then clearly $M_{\mathscr{F}}=0$ because the inductive limit of the zero system is zero. Conversely, if $M_{\mathscr{F}}=0$ then $M=\Ker(j_{M})=\mathscr{F}(M)$. $\Box$ \\

Now we prove a very useful result:\\

\begin{proposition}\label{G-lemma 2} Let $f:M\rightarrow N$ and $g:M\rightarrow P$ be $R-$linear maps such that $\Ker(g)$ and $\Coker(g)$ are $\mathscr{F}-$negligible. Then there is a unique $R-$linear map $h:P\rightarrow N_{\mathscr{F}}$ which making commutative the following diagram $$\xymatrix{
M\ar[r]^{f} \ar[d]^{g} & N\ar[d]^{j_{N}\:\:\:\:} \\P\ar[r]^{h\:\:\:\: } & N_{\mathscr{F}}.}$$ \\
\end{proposition}

{\bf Proof.} For each $x\in P$ let $J=\Ann_{R}\big(x+\Ima(g)\big)$ which belongs to $\mathscr{F}$ since $\Coker(g)$ is $\mathscr{F}-$negligible. Then consider the map $h_{x}:J\rightarrow N/\mathscr{F}(N)$ which maps each $a\in J$ into $f(m)+\mathscr{F}(N)$ where $m$ is an element of $M$ such that $g(m)=ax$. The map $h_{x}$ is well-defined since $\Ker(g)$ is $\mathscr{F}-$negligible. It is also $R-$linear. Therefore $h_{x}\in\Hom_{R}\big(J, N/\mathscr{F}(N)\big)$. We define $h:P\rightarrow N_{\mathscr{F}}$ as $x\rightsquigarrow [h_{x}]$. Clearly the map $h$ is $R-$linear and the completed diagram is commutative. Suppose $\psi:P\rightarrow N_{\mathscr{F}}$
is another $R-$linear map which making commutative the foregoing diagram.
Take $x\in P$, and let $\psi(x)=[h']$ where $h':I\rightarrow N/\mathscr{F}(N)$ is a $R-$linear map for some $I\in\mathscr{F}$. It suffices to show that $h'|_{I\cap J}=(h_{x})|_{I\cap J}$ where $J=\Ann_{R}\big(x+\Ima(g)\big)$.
For each $b\in I\cap J$ there is some $m\in M$ such that $g(m)=bx$. But $\psi\big(g(m)\big)=j_{N}\big(f(m)\big)$. Therefore there is an ideal $L\in\mathscr{F}$ contained in $I$ such that for each $c\in L$, $ch'(b)=cf(m)+\mathscr{F}(N)$. Let $h'(b)=n+\mathscr{F}(N)$ then we observe that $L\subseteq\Ann_{R}\big(f(m)-n+\mathscr{F}(N)\big)$. Thus $f(m)-n+\mathscr{F}(N)\in\mathscr{F}\big(N/\mathscr{F}(N)\big)=0$. $\Box$ \\

\begin{corollary}\label{G-coro 1} For each $R-$linear map $u:M\rightarrow N$ then there is a unique $R-$linear map $u_{\mathscr{F}}:M_{\mathscr{F}}\rightarrow N_{\mathscr{F}}$ such that the following diagram is commutative
$$\xymatrix{
M\ar[r]^{u} \ar[d]^{j_{M}} & N\ar[d]^{j_{N}\:\:\:\:} \\M_{\mathscr{F}}\ar[r]^{u_{\mathscr{F}}} & N_{\mathscr{F}}.}$$
In particular, $(j_{M})_{\mathscr{F}}=j_{M_{\mathscr{F}}}$.\\
\end{corollary}

{\bf Proof.} By Lemma \ref{G-lemma 3}, $\Ker(j_{M})$ and $\Coker(j_{M})$ are $\mathscr{F}-$negligible therefore the first part of the assertion is an immediate consequence of Proposition \ref{G-lemma 2}. Take $N=M_{\mathscr{F}}$ and $u=j_{M}$ then the second part implies. $\Box$ \\

For every $R-$linear map $u:M\rightarrow N$ and for each $[f]\in M_{\mathscr{F}}$, we have $u_{\mathscr{F}}([f])=[\overline{u}\circ f]$ where $\overline{u}:M/\mathscr{F}(M)\rightarrow N/\mathscr{F}(N)$ is induced by $u$.\\

\begin{lemma}\label{G-lemma 6} If $\xymatrix{0\ar[r]&M'\ar[r]^{u}&M\ar[r]^{v}&M''}$ is an exact sequence of $R-$modules then the sequence $\xymatrix{0\ar[r]&M'_{\mathscr{F}}\ar[r]^{u_{\mathscr{F}}}&
M_{\mathscr{F}}\ar[r]^{v_{\mathscr{F}}}&M''_{\mathscr{F}}}$ is exact. \\
\end{lemma}

{\bf Proof.} Let $[f]\in\Ker(v_{\mathscr{F}})$ where $f:I\rightarrow M/\mathscr{F}(M)$ is a $R-$linear map for some $I\in\mathscr{F}$. Thus there is an ideal $J\in\mathscr{F}$ contained in $I$ such that $\overline{v}\circ f|_{J}=0$. We have $L=J\cap f^{-1}\big(\Ima(\overline{u})\big)\in\mathscr{F}$. Because for each $b\in J$, $\Ann_{R}\big(v(m)\big)\in\mathscr{F}$ where $f(b)=m+\mathscr{F}(M)$. But $\Ann_{R}\big(v(m)\big)\subseteq\Ann_{R}\big(b+L\big)$ and so $J/L$ is $\mathscr{F}-$negligible. Hence
$L\in\mathscr{F}.\mathscr{F}=\mathscr{F}$.
Now consider the map $g:L\rightarrow M'/\mathscr{F}(M')$ given by $b\rightsquigarrow m'+\mathscr{F}(M')$ where $f(b)=\overline{u}\big(m'+\mathscr{F}(M')\big)$. The map $g$ is well-defined since $\overline{u}$ is injective. It is also $R-$linear. Thus $[g]\in M'_{\mathscr{F}}$.
Clearly $f|_{L}=\overline{u}\circ g$ and so
$[f]\in\Ima(u_{\mathscr{F}})$. $\Box$ \\

\begin{lemma}\label{G-prop 1} Let $M$ be a $R-$module. Then the canonical map $(j_{M})_{\mathscr{F}}=j_{M_{\mathscr{F}}}:M_{\mathscr{F}}
\rightarrow(M_{\mathscr{F}})_{\mathscr{F}}$ is bijective. \\
\end{lemma}

{\bf Proof.} By Corollary \ref{G-coro 1}, $(j_{M})_{\mathscr{F}}=j_{M_{\mathscr{F}}}$. By Lemma \ref{G-lemma 6}, from the exact sequence $\xymatrix{0\ar[r]&\mathscr{F}(M)\ar[r]^{i}&M\ar[r]^{j_{M}}
&M_{\mathscr{F}}}$ we obtain the following exact sequence $\xymatrix{0\ar[r]&
(\mathscr{F}(M))_{\mathscr{F}}\ar[r]^{\:\:\:\:\:\:\:\:\:i_{\mathscr{F}}}
&M_{\mathscr{F}}\ar[r]^{(j_{M})_{\mathscr{F}}\:\:\:\:}
&(M_{\mathscr{F}})_{\mathscr{F}}}$. By Lemma \ref{G-lemma 7}, $\big(\mathscr{F}(M)\big)_{\mathscr{F}}=0$ hence $j_{M_{\mathscr{F}}}$ is injective thus $\mathscr{F}(M_{\mathscr{F}})=0$. Therefore $(M_{\mathscr{F}})_{\mathscr{F}}=
\colimit_{I\in\mathscr{F}}\Hom_{R}(I,M_\mathscr{F})$. Take $[f]\in(M_{\mathscr{F}})_{\mathscr{F}}$ where $f:I\rightarrow M_\mathscr{F}$ is a $R-$linear map and $I\in\mathscr{F}$. We claim that $L=f^{-1}\big(\Ima j_{M}\big)\in\mathscr{F}$. For each $a\in I$ there is an ideal $J\in\mathscr{F}$ and also there is a $R-$linear map $h:J\rightarrow M/\mathscr{F}(M)$ such that $f(a)=[h]$. To prove the claim, first we show that $J\subseteq\Ann_{R}(a+L)$. Let $b\in J$ and let $h(b)=m+\mathscr{F}(M)$ where $m\in M$. We have $f(ab)=bf(a)=[b.h]$. But $b.h=(\delta_{\overline{m}})|_{J}$ therefore $f(ab)=[\delta_{\overline{m}}]\in\Ima j_{M}$. This implies that $I/L$ is $\mathscr{F}-$negligible and so $L\in\mathscr{F}.\mathscr{F}=\mathscr{F}$. This establishes the claim.
Now consider the map $g:L\rightarrow M/\mathscr{F}(M)$ which maps each $a\in L$ into $m+\mathscr{F}(M)$ where $f(a)=j_{M}(m)$. The map $g$ is well-defined since $\Ker(j_{M})=\mathscr{F}(M)$.
It is also $R-$linear. Hence $[g]\in M_{\mathscr{F}}$. Clearly $f|_{L}=\overline{j_{M}}\circ g$ thus $[f]=j_{M_{\mathscr{F}}}([g])$ and so $j_{M_{\mathscr{F}}}$ is surjective. $\Box$ \\

Now we are ready to prove the first main result of this article:\\

\begin{theorem}\label{G-lemma 10} Let $\phi:M\rightarrow N$ be a $R-$linear map such that $\Ker\phi$ and $\Coker\phi$ are $\mathscr{F}-$negligible. Then $\phi_{\mathscr{F}}:M_{\mathscr{F}}\rightarrow N_{\mathscr{F}}$ is bijective.\\
\end{theorem}

{\bf Proof.} By Proposition \ref{G-lemma 2}, there is a (unique) $R-$linear map $h:N\rightarrow M_{\mathscr{F}}$ such that $h\circ\phi=j_{M}$. By Lemma \ref{G-prop 1}, $j_{M_{\mathscr{F}}}$ is bijective. We shall prove that
$j^{-1}_{M_{\mathscr{F}}}\circ h_{\mathscr{F}}$ is the inverse of $\phi_{\mathscr{F}}$. We have $j^{-1}_{M_{\mathscr{F}}}\circ h_{\mathscr{F}}\circ\phi_{\mathscr{F}}=j^{-1}_{M_{\mathscr{F}}}
\circ(h\circ\phi)_{\mathscr{F}}=\Identity_{M_{\mathscr{F}}}$.
To conclude the proof, by Corollary \ref{G-coro 1}, it suffices to show that $\phi_{\mathscr{F}}\circ j^{-1}_{M_{\mathscr{F}}}\circ h_{\mathscr{F}}\circ j_{N}=j_{N}$. But $\phi_{\mathscr{F}}\circ j^{-1}_{M_{\mathscr{F}}}\circ h_{\mathscr{F}}\circ j_{N}=\phi_{\mathscr{F}}\circ j^{-1}_{M_{\mathscr{F}}}\circ j_{M_{\mathscr{F}}}\circ h=\phi_{\mathscr{F}}\circ h$. For each $x\in N$, from the proof of Proposition \ref{G-lemma 2}, we know that $h(x)=[h_{x}]$ where $h_{x}:J=\Ann_{R}(x+\Ima\phi)\rightarrow M/{\mathscr{F}}(M)$ is a $R-$linear map which maps each $a\in J$ into $m+\mathscr{F}(M)$
such that $\phi(m)=ax$. Thus $\overline{\phi}\circ h_{x}=(\delta_{\overline{x}})|_{J}$ where $\overline{x}=x+\mathscr{F}(N)$. This means that $\phi_{\mathscr{F}}\circ h=j_{N}$. $\Box$ \\

Now the G-localization rings are introduced. Let $M$ be a $R-$module. For each $[f]\in R_{\mathscr{F}}$ and for each $[g]\in M_{\mathscr{F}}$ where $f:I\rightarrow R/\mathscr{F}(R)$ and $g:J\rightarrow M/\mathscr{F}(M)$ are $R-$linear maps and $I,J\in\mathscr{F}$, we have $f^{-1}(\overline{J})\in\mathscr{F}$ where $\overline{J}=J+\mathscr{F}(R)/\mathscr{F}(R)$. Moreover $g$ induces a map $\overline{J}\rightarrow R/\mathscr{F}(R)$ given by $b+\mathscr{F}(R)\rightsquigarrow g(b)$ which we denote it by $\overline{g}$. The map $\overline{g}$ is clearly well-defined. Now we define the pairing $R_{\mathscr{F}}\times M_{\mathscr{F}}\rightarrow M_{\mathscr{F}}$ as  $[f].[g]=[\overline{g}\circ\big(f|_{f^{-1}(\overline{J})}\big)]$. It is easy to see that it is $R-$bilinear (details omitted). \\

In particular, the binary operation $R_{\mathscr{F}}\times R_{\mathscr{F}}\rightarrow R_{\mathscr{F}}$, as multiplicative, puts a commutative ring structure on the $R-$module $R_{\mathscr{F}}$. Its commutativity implies from the fact that for every $R-$linear maps $f,g:I\rightarrow R/\mathscr{F}(R)$ where $I$ is an ideal of $R$ and for every elements $a,b\in I$ then $\overline{g}\big(f(ab)\big)=\overline{f}\big(g(ab)\big)$. The unit element of the ring $R_{\mathscr{F}}$ is $[\pi]$ where $\pi:R\rightarrow R/\mathscr{F}(R)$ is the canonical map. The map $j_{R}:R\rightarrow R_{\mathscr{F}}$ is a ring homomorphism. Moreover, the pairing
$R_{\mathscr{F}}\times M_{\mathscr{F}}\rightarrow M_{\mathscr{F}}$ puts a $R_{\mathscr{F}}-$module structure on $M_{\mathscr{F}}$. For every $R-$linear map $u:M\rightarrow N$ then $u_{\mathscr{F}}:M_{\mathscr{F}}\rightarrow N_{\mathscr{F}}$ is $R_{\mathscr{F}}-$linear. In fact $M\rightsquigarrow M_{\mathscr{F}}$ is a left exact functor from the category of $R-$modules into the category of $R_{\mathscr{F}}-$modules, see Corollary \ref{G-coro 1} and Lemma \ref{G-lemma 6}. It is called the Gabriel localization (G-localization) functor with respect to the system $\mathscr{F}$.\\

\begin{proposition}\label{lemma 189} Let $J$ be an ideal of $R$ such that $JR_{\mathscr{F}}=R_{\mathscr{F}}$. Then $J\in\mathscr{F}$.\\
\end{proposition}

{\bf Proof.} We may write $1=\sum\limits_{i=1}^{n}z_{i}j_{R}(a_{i})$ where $a_{i}\in J$ and $z_{i}\in R_{\mathscr{F}}$ for all $i$. Let $f:I\rightarrow N/\mathscr{F}(N)$ be a $R-$linear map where $N=R/J$ and $I\in\mathscr{F}$. Then we have $[f]=\sum\limits_{i=1}^{n}z_{i}.[a_{i}.f]=0$ since $a.f=0$ for all $a\in J$. Therefore $R/J$ is $\mathscr{F}-$negligible and so $J\in\mathscr{F}.\mathscr{F}=\mathscr{F}$. $\Box$ \\

Note that the converse of Proposition \ref{lemma 189} does not necessarily hold. In fact, the condition ``$IR_{\mathscr{F}}=R_{\mathscr{F}}$ for all $I\in\mathscr{F}$'', as we shall observe in the article, is a crucial point of the G-localization rings. Many interesting and major facts are equivalent or imply from this condition, see for example, Corollary \ref{coro 6533} and Theorem \ref{th621890}.\\

The second main result of this article is the following.\\

\begin{theorem}\label{G-th 5} Let $\mathscr{F}$ be an idempotent topologizing system on the ring $R$ and let $\phi:R\rightarrow S$ be a ring map such that $IS=S$ for all $I\in\mathscr{F}$. Then the following conditions hold.\\
$\textbf{(i)}$ For each $S-$module $M$, the canonical map $j_{M}:M\rightarrow M_{\mathscr{F}}$ is bijective.\\
$\textbf{(ii)}$ The map $\psi=j^{-1}_{S}\circ\phi_{\mathscr{F}}$
is the only ring homomorphism from $R_{\mathscr{F}}$ into $S$ such that $\phi=\psi\circ j_{R}$.\\
\end{theorem}

{\bf Proof.} $\textbf{(i)}:$ If $m\in\mathscr{F}(M)$ then $I=\Ann_{R}(m)\in\mathscr{F}$. Thus we may write $1=\sum\limits_{i}s'_{i}\phi(a_{i})$ where $a_{i}\in I$
and $s'_{i}\in S$. Therefore $m=\sum\limits_{i}s'_{i}\phi(a_{i})m=0$ hence $j_{S}$ is injective. Let $f:J\rightarrow M$ be a $R-$linear map where $J\in\mathscr{F}$. We may write $1=\sum\limits_{j=1}^{n}s_{j}\phi(b_{j})$ where $b_{j}\in J$ and $s_{j}\in S$ for all $j$. For each $c\in J$ we have $f(c)=
\sum\limits_{j=1}^{n}s_{j}\phi(b_{j})f(c)=\sum\limits_{j=1}^{n}s_{j}f(b_{j}c)=
\big(\sum\limits_{j=1}^{n}s_{j}f(b_{j})\big)\phi(c)$. Therefore $(\delta_{m})|_{J}=f$ where
$m=\sum\limits_{j=1}^{n}s_{j}f(b_{j})$ and $\delta_{m}:R\rightarrow M$ which maps each $r\in R$ into $r.m=\phi(r)m$. This means that $j_{M}(m)=[f]$ and so $j_{M}$ is surjective. \\
$\textbf{(ii)}:$ First we show that the map $\psi=j_{S}^{-1}\circ\phi_{\mathscr{F}}$ is actually a ring homomorphism. Clearly it is additive and transforms the unit element of $R_{\mathscr{F}}$ to the unit of $S$. To prove that it is multiplicative take two elements $z,z'\in R_{\mathscr{F}}$ with representations $f$ and $g$, i.e., $z=[f]$ and $z'=[g]$. By passing to the restrictions, if it is necessary, therefore we may assume that the $R-$linear maps $f$ and $g$ have the same domain. Namely $f,g:J\rightarrow R/\mathscr{F}(R)$ where $J\in\mathscr{F}$. We know that $f^{-1}(\overline{J})\in\mathscr{F}$ where $\overline{J}=J+\mathscr{F}(R)/\mathscr{F}(R)$.
Therefore we may write $1=\sum\limits_{j}\phi(b_{j})s_{j}$ where $b_{j}\in f^{-1}(\overline{J})\subseteq J$ and $s_{j}\in S$. We have then $\psi(z.z')=
\sum\limits_{j}\phi(c'_{j})s_{j}$ where for each $j$, $f(b_{j})=c_{j}+\mathscr{F}(R)$, $g(c_{j})=c'_{j}+\mathscr{F}(R)$ and $c_{j}\in J$. Similarly $\psi(z)=\sum\limits_{j}\phi(c_{j})s_{j}$ and $\psi(z')=\sum\limits_{i}\phi(e_{i})s_{i}$ where for each $i$, $g(b_{i})=e_{i}+\mathscr{F}(R)$. Thus $\psi(z)
\psi(z')=
\sum\limits_{j}\big(\sum\limits_{i}\phi(e_{i}c_{j})s_{i}\big)s_{j}$. But for each $j$ we have $\phi(c'_{j})=\sum\limits_{i}\phi(c'_{j}b_{i})s_{i}=
\sum\limits_{i}\overline{\phi}\big(c'_{j}b_{i}+\mathscr{F}(R)\big)s_{i}=
\sum\limits_{i}\overline{\phi}\big(b_{i}g(c_{j})\big)s_{i}=
\sum\limits_{i}\overline{\phi}\big(c_{j}g(b_{i})\big)s_{i}=
\sum\limits_{i}\overline{\phi}\big(e_{i}c_{j}+\mathscr{F}(R)\big)s_{i}=
\sum\limits_{i}\phi(e_{i}c_{j})s_{i}$. Therefore $\psi(z.z')=\psi(z)
\psi(z')$. Finally, suppose $\psi':R_{\mathscr{F}}\rightarrow S$ is another ring map such that $\psi'\circ j_{R}=\phi$. Clearly $j_{S}\circ\psi'$ is $R-$linear and $(j_{S}\circ\psi')\circ j_{R}=j_{S}\circ\phi$. Therefore, by Corollary \ref{G-coro 1}, $j_{S}\circ\psi'=\phi_{\mathscr{F}}$ and so $\psi'=j_{S}^{-1}\circ\phi_{\mathscr{F}}$.  $\Box$ \\

\section{Flat epimorphisms as G-localizations}

By an epimorphism $\phi:R\rightarrow S$ we mean it is an epimorphism in the category of commutative rings.  We should mention that the surjective ring maps are just special cases of the epimorphisms. For example, the canonical ring map $\mathbb{Z}\rightarrow\mathbb{Q}$ is an epimorphism while it is not surjective. A ring map which is both flat and an epimorphism is called a flat epimorphism. The canonical map $R\rightarrow S^{-1}R$ where $S$ is a multiplicative subset of a ring $R$ is a typical example of  flat epimorphisms.\\

The following lemma is a well-known result but we have provided a proof for the sake of completeness. \\

\begin{lemma}\label{G-lemma 9} Let $\phi:R\rightarrow S$ be a ring homomorphism. Then the following conditions are equivalent.\\
$\mathbf{(i)}$ $\phi$ is an epimorphism.\\
$\mathbf{(ii)}$ For each $s\in S$, $s\otimes1=1\otimes s$.\\
$\mathbf{(iii)}$ The map $p: S\otimes_{R}S\rightarrow S$ defined by $s\otimes s'\rightsquigarrow ss'$ is bijective.\\
$\mathbf{(iv)}$ The map $j: S\rightarrow S\otimes_{R}S$ defined by $s\rightsquigarrow 1\otimes s$ is bijective.\\
\end{lemma}

{\bf Proof.} $\mathbf{(i)}\Rightarrow\mathbf{(ii):}$ We have  $i\circ\phi=j\circ\phi$ where $i,j:S\rightarrow S\otimes_{R}S$ are the canonical ring maps which map each $s\in S$ into $s\otimes1$ and $1\otimes s$, respectively. \\
$\mathbf{(ii)}\Rightarrow\mathbf{(iii):}$ We have $\Ker(p)=\langle s\otimes1-1\otimes s : s\in S\rangle$ because if $ss'=0$ then we may write $s\otimes s'=1\otimes s'(s\otimes1- 1\otimes s)$.\\
$\mathbf{(iii)}\Rightarrow\mathbf{(i):}$ Let $f,g:S\rightarrow T$ be ring maps such that $f\circ\phi=g\circ\phi$. By the universal property of the pushouts, there is a (unique) ring map $\psi:S\otimes_{R}S\rightarrow T$ such that $f=\psi\circ i$ and $g=\psi\circ j$. But $i=j$ since $\Ker(p)=0$. Thus $f=g$.\\ $\mathbf{(ii)}\Rightarrow\mathbf{(iv):}$ We have $s\otimes s'=(s\otimes1)(1\otimes s')=1\otimes ss'=j(ss')$.\\
$\mathbf{(iv)}\Rightarrow\mathbf{(iii):}$ We have $p\circ j=\Identity$. Therefore $p$ is bijective. $\Box$ \\

The third main result of this article is the following.\\

\begin{theorem}\label{G-lemma 1} Let $\phi:R\rightarrow S$ be a ring map and let $\mathscr{F}$ be the set of ideals of $R$ whose extensions under $\phi$ are equal to $S$. Then the following conditions hold.\\
$\textbf{(i)}$ The family $\mathscr{F}$ is an idempotent toplogizing system on the ring $R$. \\
$\textbf{(ii)}$ If $\phi$ is a flat epimorphism then for each $R-$module $M$, the map $\eta_{\mathscr{F}}:M_{\mathscr{F}}
\rightarrow(S\otimes_{R}M)_{\mathscr{F}}$ induced by the canonical map $\eta:M\rightarrow S\otimes_{R}M$ given by $m\rightsquigarrow1\otimes m$ is an isomorphism. In particular the map $\phi_{\mathscr{F}}$ is bijective.\\
\end{theorem}

{\bf Proof.} $\textbf{(i)}:$ Clearly $\mathscr{F}$ is non-empty since $R\in\mathscr{F}$. Suppose $I\in\mathscr{F}$ and $J$ is an ideal of $R$ such that $J:a\in\mathscr{F}$ for all $a\in I$. We may write $1=\sum\limits_{j=1}^{n}s_{j}\phi(a_{j})$
where $a_{j}\in I$ and $s_{j}\in S$. For each $j$, we may also write $1=\sum\limits_{i_{j}=1}^{k_{j}}s_{i_{j},j}\phi(b_{i_{j},j})$ where $b_{i_{j},j}\in J:a_{j}$ and $s_{i_{j},j}\in S$. We have then $1=\sum\limits_{i_{1}=1}^{k_{1}}...\sum\limits_{i_{n}=1}^{k_{n}}
\sum\limits_{j=1}^{n}s_{i_{1},1}...s_{i_{n},n}s_{j}
\phi(b_{i_{1},1}...b_{i_{n},n}a_{j})$ and all the elements $b_{i_{1},1}...b_{i_{n},n}a_{j}$ belong to $J$ and so $J\in\mathscr{F}$. Therefore, by Lemma \ref{G-lemma 5}, $\mathscr{F}$ is an idempotent topologizing system. \\
$\textbf{(ii)}:$ To prove the assertion, by Theorem \ref{G-lemma 10}, it suffices to show that $\Ker\eta$ and $\Coker\eta$ are $\mathscr{F}-$negligible. Let $m\in\Ker\eta$ then
the map $\lambda:S\otimes_{R}N\rightarrow S\otimes_{R}M$ induced by the canonical injection $N=Rm\rightarrow M$ is injective since $S$ is $R-$flat. But $\Ima\lambda=0$ because $1\otimes m=0$. This implies that $\Ann_{R}(m)S=S$ and so $\Ann_{R}(m)\in\mathscr{F}$. Thus $\Ker\eta$ is $\mathscr{F}-$negligible. \\ By applying the right exact functor $S\otimes_{R}-$ to the exact sequence $\xymatrix{M\ar[r]^{\eta\:\:\:\:\:\:\:\:}&S\otimes_{R}
M\ar[r]&\Coker\eta\ar[r]&0}$ we then obtain the following exact sequence
$\xymatrix{S\otimes_{R}M\ar[r]^{1\otimes\eta\:\:\:\:\:\:\:\:\:\:\:\:}&S\otimes_{R}
(S\otimes_{R}M)\ar[r]&S\otimes_{R}\Coker\eta\ar[r]&0.}$ \\
The map $1\otimes\eta$ factors as $\xymatrix{S\otimes_{R}M\ar[r]^{j\otimes 1\:\:\:\:\:\:\:\:\:\:\:\:}&
(S\otimes_{R}S)\otimes_{R}M\ar[r]^{\:\:\:\:\:\:\:\:\:\:\:\:\:\:\:\:\:\:\:\:\simeq}&}$\\
$\xymatrix{S\otimes_{R}(S\otimes_{R}M).}$
By Lemma \ref{G-lemma 9}, $j\otimes1$ is bijective, hence so is
$1\otimes\eta$. In particular it is surjective and so $S\otimes_{R}\Coker\eta=0$. For each $x\in\Coker\eta$, the map
$S\otimes_{R}Rx\rightarrow S\otimes_{R}\Coker\eta$ induced by the canonical injection $Rx\rightarrow\Coker\eta$ is injective
since $S$ is $R-$flat. This implies that $\Ann_{R}(x)S=S$ and so $\Ann_{R}(x)\in\mathscr{F}$. Therefore $\Coker\eta$ is also $\mathscr{F}-$negligible. Finally, the map $\phi$ factors as $\xymatrix{R\ar[r]^{\eta\:\:\:\:\:\:\:}&R\otimes_{R}S\ar[r]^{\:\:\:\:\:\simeq}&S}$
hence $\phi_{\mathscr{F}}$ is bijective. $\Box$ \\

\begin{corollary}\label{G-coro 2} Let $\phi:R\rightarrow S$ be a ring map and let $\mathscr{F}$ be the set of ideals $I$ of $R$ such that $IS=S$. If $\phi$ is a flat epimorphism then there is a unique isomorphism of rings $\psi:R_{\mathscr{F}}\rightarrow S$ such that $\phi=\psi\circ j_{R}$. \\
\end{corollary}

{\bf Proof.} By Theorem \ref{G-th 5}, the map $\psi=j^{-1}_{S}\circ\phi_{\mathscr{F}}$ is the only ring homomorphism such that $\phi=\psi\circ j_{R}$. By Theorem \ref{G-lemma 1},
$\phi_{\mathscr{F}}$ is bijective, therefore $\psi$ is an isomorphism. $\Box$ \\

\begin{corollary}\label{coro 87} Let $J$ be an ideal of a ring $R$. Then $R/J$ is $R-$flat if and only if $\Ann_{R}(a)+J=R$ for all $a\in J$.\\
\end{corollary}

{\bf Proof.} Suppose $R/J$ is $R-$flat. Let $\mathscr{F}$ be the set of ideals $I$ of $R$ such that $I+J=R$. By Corollary \ref{G-coro 2}, there is a (unique) isomorphism of rings $\psi: R/J\rightarrow R_{\mathscr{F}}$ such that $j_{R}=\psi\circ\pi$ where $\pi:R\rightarrow R/J$ is the canonical map. We have $\mathscr{F}(R)=\Ker(j_{R})=j_{R}^{-1}(0)=\pi^{-1}\big(\psi^{-1}(0)\big)
=\Ker(\pi)=J$. Conversely, let $f:M\rightarrow N$ be an injective $R-$linear map. To prove the assertion it suffices to show that the induced map $M/JM\rightarrow N/JN$ given by $m+JM\rightsquigarrow f(m)+JN$ is injective. If $f(m)\in JN$ then we may write $f(m)=\sum\limits_{i=1}^{s}a_{i}n_{i}$ where $a_{i}\in J$ and $n_{i}\in N$ for all $i$. By the hypothesis, there are elements $b_{i}\in\Ann_{R}(a_{i})$ and $c_{i}\in J$ such that $1=b_{i}+c_{i}$. It follows that $1=(b_{1}+c_{1})(b_{2}+c_{2})...(b_{s}+c_{s})=b+c$ where $b=b_{1}b_{2}...b_{s}$ and $c\in J$. Thus $f(m)=bf(m)+cf(m)=f(cm)$. Therefore $m=cm\in JM$. $\Box$ \\

Corollary \ref{coro 87}, in particular, tells us that if the ideal $J$ has a generating set $S$ such that each $a\in S$ can be written as $a=a^{2}b$ for some $b\in R$ then $R/J$ is $R-$flat. As another application, if $J$ is an ideal of a domain $R$ such that $R/J$ is $R-$flat then we have either $J=0$ or $J=R$.\\

\begin{lemma}\label{coro 652} Let $\mathscr{F}$ be an idempotent topologizing system on the ring $R$ such that $IR_{\mathscr{F}}=R_{\mathscr{F}}$ for all $I\in\mathscr{F}$. Then $j_{R}$ is an epimorphism.\\
\end{lemma}

{\bf Proof.} Let $f,g:R_{\mathscr{F}}\rightarrow S$ be two ring maps such that $f\circ j_{R}=g\circ j_{R}$. If we consider the ring map $\phi=f\circ j_{R}:R\rightarrow S$ then $IS=S$ for all $I\in\mathscr{F}$. Therefore, by Theorem \ref{G-th 5}, there is a unique ring map $\psi:R_{\mathscr{F}}\rightarrow S$ such that $\phi=\psi\circ j_{R}$. Thus $f=g$. $\Box$ \\

The converse of Corollary \ref{G-coro 2} also holds even under a mild hypothesis which is another main result of this article:\\

\begin{theorem}\label{remark 200} Let $\mathscr{F}$ be an idempotent topologizing system on the ring $R$ and let $\phi:R\rightarrow S$ be a ring map such that $IS=S$ for all $I\in\mathscr{F}$. If there is an isomorphism of rings $\psi:R_{\mathscr{F}}\rightarrow S$ such that $\phi=\psi\circ j_{R}$ then $\phi:R\rightarrow S$ is a flat epimorphism and $\psi=j_{S}^{-1}\circ\phi_{\mathscr{F}}$.\\
\end{theorem}

{\bf Proof.} We have $IR_{\mathscr{F}}=R_{\mathscr{F}}$ for all $I\in\mathscr{F}$. Therefore, by Lemma \ref{coro 652}, $j_{R}$ and so $\phi$ are epimorphisms. By Theorem \ref{G-th 5}, $\psi=j_{S}^{-1}\circ\phi_{\mathscr{F}}$. Therefore $\phi_{\mathscr{F}}$ is bijective. This, in particular, implies that $j_{S}(s)\in\Ima\phi_{\mathscr{F}}$ for all $s\in S$. Thus there is an ideal $I\in\mathscr{F}$ such that $I\subseteq\Ann_{R}(s+\Ima\phi)$. Therefore $\Coker\phi$ and so any finite direct sum of it are $\mathscr{F}-$negligible. To prove that $S$ is $R-$flat, by \cite[Theorem 7.7]{Matsumura}, it suffices to show that for each ideal $J$ of $R$ then the canonical map $J\otimes_{R}S\rightarrow S$ given by $\sum\limits_{i=1}^{n}a_{i}\otimes s_{i}\rightsquigarrow\sum\limits_{i=1}^{n}\phi(a_{i})s_{i}$ is injective.
Suppose $\sum\limits_{i=1}^{n}\phi(a_{i})s_{i}=0$ where $a_{i}\in J$ and $s_{i}\in S$ for all $i$. We have $\Ann_{R}(x)\in\mathscr{F}$ where $x=(s_{i}+\Ima\phi)^{n}_{i=1}\in(\Coker\phi)^{n}$ since $(\Coker\phi)^{n}$ is $\mathscr{F}-$negligible. Thus there are elements $b_{1},...,b_{m}\in\Ann_{R}(x)$ and also elements $s'_{1},...,s'_{m}\in S$ such that $1=\sum\limits_{j=1}^{m}\phi(b_{j})s'_{j}$. Moreover there are elements $r_{i,j}\in R$ such that $s_{i}\phi(b_{j})=\phi(r_{i,j})$ for all $i,j$.
Thus $c_{j}=\sum\limits_{i=1}^{n}a_{i}r_{i,j}\in\Ker\phi=\mathscr{F}(R)$ for all $j$. Hence for each $j$ there are elements $r'_{j,1},...,r'_{j,N}\in\Ann_{R}(c_{j})$ and also elements $s''_{j,1},...,s''_{j,N}\in S$ such that $1=\sum\limits_{k=1}^{N}\phi(r'_{j,k})s''_{j,k}$.
Now we have $\sum\limits_{i=1}^{n}a_{i}\otimes s_{i}=\sum\limits_{i=1}^{n}a_{i}\otimes\big(
\sum\limits_{j=1}^{m}s_{i}\phi(b_{j})s'_{j}\big)=
\sum\limits_{j=1}^{m}c_{j}\otimes s'_{j}$. For each $j$,
$c_{j}\otimes s'_{j}=c_{j}\otimes
\big(\sum\limits_{k=1}^{N}\phi(r'_{j,k})s'_{j}s''_{j,k}\big)=
\sum\limits_{k=1}^{N}r'_{j,k}c_{j}\otimes s'_{j}s''_{j,k}=0$. $\Box$ \\

\begin{corollary}\label{coro 6533} Let $\mathscr{F}$ be an idempotent topologizing system on the ring $R$ such that $IR_{\mathscr{F}}=R_{\mathscr{F}}$ for all $I\in\mathscr{F}$. Then $j_{R}$ is a flat epimorphism.\\
\end{corollary}

{\bf Proof.} It is an immediate consequence of Theorem \ref{remark 200}. $\Box$ \\

It is natural to ask whether the converse of Corollary \ref{coro 6533} holds.\\

\section{G-localizations of finite type systems}

\begin{definition} An $R-$module $M$ is said to be $\mathscr{F}-$closed if the canonical map $j_{M}$ is bijective. It is called strongly $\mathscr{F}-$closed if the canonical map
$M\rightarrow\Hom_{R}(I,M)$ given by $m\rightsquigarrow(\delta_{m})|_{I}$ is bijective for all $I\in\mathscr{F}$. \\
\end{definition}

For each $R-$module $M$ then $M_{\mathscr{F}}$, by Lemma \ref{G-prop 1}, is $\mathscr{F}-$closed. Clearly each strongly $\mathscr{F}-$closed module is $\mathscr{F}-$closed. By the category of  $\mathscr{F}-$closed modules we mean a full subcategory of the category of $R-$modules whose objects are the $\mathscr{F}-$closed modules. \\

\begin{lemma}\label{prop 431} Let $\mathscr{F}$ be an idempotent topologizing system on the ring $R$ such that $IR_{\mathscr{F}}=R_{\mathscr{F}}$ for all $I\in\mathscr{F}$. Then every $\mathscr{F}-$closed module is strongly $\mathscr{F}-$closed and the G-localization functor with respect to $\mathscr{F}$ is an equivalence between the category of $\mathscr{F}-$closed modules and the category of $R_{\mathscr{F}}-$modules. \\
\end{lemma}

{\bf Proof.} Let $N$ be a $\mathscr{F}-$closed module and let $I\in\mathscr{F}$. Then the canonical map $N\rightarrow\Hom_{R}(I,N)$ is injective since $\mathscr{F}(N)=0$. Let $f:I\rightarrow N$ be a $R-$linear map. We may write $1=\sum\limits_{i=1}^{n}j_{R}(a_{i})z_{i}$ where $a_{i}\in I$ and $z_{i}\in R_{\mathscr{F}}$ for all $i$. There exists an element $x\in N$ such that $j_{N}(x)=\sum\limits_{i=1}^{n}z_{i}j_{N}\big(f(a_{i})\big)$. Now, for each $a\in I$, we have $j_{N}(ax)=j_{N}\big(f(a)\big)$ and so $ax=f(a)$. Thus $(\delta_{x})|_{I}=f$. Therefore $N$ is strongly $\mathscr{F}-$closed. Moreover, for each $R-$module $M$ then the map $\Hom_{R}(M,N)\rightarrow\Hom_{R_{\mathscr{F}}}
(M_{\mathscr{F}},N_{\mathscr{F}})$ given by $u\rightsquigarrow u_{\mathscr{F}}$ is bijective. Because suppose $u_{\mathscr{F}}=u'_{\mathscr{F}}$ where $u,u':M\rightarrow N$ are $R-$linear maps. We have $u=j_{N}^{-1}\circ u_{\mathscr{F}}\circ j_{M}=u'$.
Let $\phi:M_{\mathscr{F}}\rightarrow N_{\mathscr{F}}$ be a $R-$linear map. Then, by Corollary \ref{G-coro 1}, $u_{\mathscr{F}}=\phi$ where $u=j_{N}^{-1}\circ\phi\circ j_{M}$. Finally, we show that the G-localization functor is essentially surjective. Let $L$ be a $R_{\mathscr{F}}-$module. By Theorem \ref{G-th 5}, $j_{L}$ is bijective. Thus $L$ as $R-$module is $\mathscr{F}-$closed. It remains to show that $j_{L}$ is $R_{\mathscr{F}}-$linear. We have $j_{L}=\sigma_{L}\circ\eta$ where $\eta:L\rightarrow R_{\mathscr{F}}\otimes_{R}L$ is the canonical map. The map $\sigma_{L}$ is already $R_{\mathscr{F}}-$linear. By Lemma \ref{coro 652}, $j_{R}$ is an epimorphism. Note that for any epimorphism of rings $\phi:R\rightarrow S$ and for any $S-$modules $M$ and $N$ then the two $S-$module structures on $M\otimes_{R}N$ defined on pure tensors by $s.(m\otimes n)=sm\otimes n$ and $s\ast(m\otimes n)=m\otimes sn$ are the same since $s\otimes1=1\otimes s$ for all $s\in S$ and so the canonical map $\eta:N\rightarrow S\otimes_{R}N$ is $S-$linear. Therefore $j_{L}$ is $R_{\mathscr{F}}-$linear. $\Box$ \\

\begin{definition} An idempotent topologizing system is called \emph{of  finite type} if every element of it containing a finitely generated ideal belonging to the system. For example, if $\phi:R\rightarrow S$ is a ring map then the system $\{I\subseteq R : IS=S\}$ is of finite type.\\
\end{definition}

In the next result we shall use, for each $R-$module $M$, the canonical $R_{\mathscr{F}}-$linear map $\sigma_{M}:R_{\mathscr{F}}\otimes_{R}M\rightarrow M_{\mathscr{F}}$ which maps each pure tensor $a\otimes m$ into $a.j_{M}(m)$ for all $a\in R_{\mathscr{F}}$ and all $m\in M$. Note that $\sigma_{R}$ is bijective and
for each $R-$linear map $\phi:M\rightarrow N$ then the following diagram is commutative $$\xymatrix{
 R_{\mathscr{F}}\otimes_{R}M\ar[r]^{1\otimes\phi} \ar[d]^{\sigma_{M}} &  R_{\mathscr{F}}\otimes_{R}N\ar[d]^{\sigma_{N}} \\ M_{\mathscr{F}}\ar[r]^{\:\:\:\:\phi_{\mathscr{F}}} & N_{\mathscr{F}}.}$$ \\

The following is another main result of this article:\\

\begin{theorem}\label{th621890} Let $\mathscr{F}$ be an idempotent topologizing system on the ring $R$. Then the following conditions are equivalent.\\
$\textbf{(i)}$ For each $R-$module $M$, the canonical map $\sigma_{M}:R_{\mathscr{F}}\otimes_{R}M\rightarrow M_{\mathscr{F}}$ is bijective.\\
$\textbf{(ii)}$ For every $R-$module $M$, $\mathscr{F}(M)$ is the kernel of the canonical map $\eta:M\rightarrow R_{\mathscr{F}}\otimes_{R}M$.\\
$\textbf{(iii)}$ For each $I\in\mathscr{F}$, $IR_{\mathscr{F}}=R_{\mathscr{F}}$. \\
$\textbf{(iv)}$ The G-localization functor with respect to $\mathscr{F}$ is exact and preserves direct sums.\\
$\textbf{(v)}$ The G-localization functor with respect to $\mathscr{F}$ is exact  and the system $\mathscr{F}$ is of finite type.\\
$\textbf{(vi)}$ The G-localization functor with respect to $\mathscr{F}$ is essentially surjective.\\
$\textbf{(vii)}$ For each $R_{\mathscr{F}}-$module $N$, $\mathscr{F}(N)=0$.\\
\end{theorem}

{\bf Proof.} $\textbf{(i)}\Rightarrow\textbf{(ii)}:$ We have $j_{M}=\sigma_{M}\circ\eta$. Thus $\mathscr{F}(M)=\Ker j_{M}=j_{M}^{-1}(0)=\eta^{-1}\big(\sigma_{M}^{-1}(0)\big)=\Ker\eta$.\\
$\textbf{(ii)}\Rightarrow\textbf{(iii)}:$ If $I\in\mathscr{F}$ then $R/I$ is $\mathscr{F}-$negligible. Thus $\Ker\eta=R/I$ where $\eta:R/I\rightarrow R_{\mathscr{F}}\otimes_{R}R/I$ is the canonical map. Therefore $\Ima\eta=0$.
This implies that $R_{\mathscr{F}}\otimes_{R}R/I=0$ since $z\otimes(r+I)=\big(z\otimes(1+I)\big)\big(1\otimes(r+I)\big)=0$. It follows that $IR_{\mathscr{F}}=R_{\mathscr{F}}$. \\
$\textbf{(iii)}\Rightarrow\textbf{(i)}:$ By Corollary \ref{coro 6533}, $j_{R}$ is a flat epimorphism. Moreover, by Proposition \ref{lemma 189}, $\mathscr{F}=\{I\subseteq R : IR_{\mathscr{F}}=R_{\mathscr{F}}\}$.
Therefore, by the proof of Theorem \ref{G-lemma 1}, $\Ker\eta$ and $\Coker\eta$ are $\mathscr{F}-$negligible where $\eta:M\rightarrow R_{\mathscr{F}}\otimes_{R}M$ is the canonical map. Thus, by Proposition \ref{G-lemma 2}, there is a unique $R-$linear map $h: R_{\mathscr{F}}\otimes_{R}M\rightarrow M_{\mathscr{F}}$ such that $j_{M}=h\circ\eta$. By Theorem \ref{G-th 5}, $j_{N}$ is bijective where $N=R_{\mathscr{F}}\otimes_{R}M$. Moreover, by Theorem \ref{G-lemma 1}, $\eta_{\mathscr{F}}$ is bijective. We also have $j_{M}=(\eta_{\mathscr{F}}^{-1}\circ j_{N})\circ\eta$. Therefore $\sigma_{M}=\eta_{\mathscr{F}}^{-1}\circ j_{N}$ and so $\sigma_{M}$ is bijective.\\
$\textbf{(i)}\Rightarrow\textbf{(iv)}$ and $\textbf{(v)}:$ Easy.\\
$\textbf{(iv)}\Rightarrow\textbf{(i)}:$ Let $L$ be a free $R-$module and let $B=\{x_{i} : i\in I\}$ be a basis of $L$. Then there is a bijective map $\psi:L\rightarrow\bigoplus\limits_{i\in I}R$
which maps each $x\in L$ into $(r_{i})_{i\in I}$ where $x=\sum\limits_{i\in I}r_{i}x_{i}$. We denote the converse of $\psi$ by $\theta$ and for a given $R-$module $M$, $\bigoplus\limits_{i\in I}M$ is denoted by $M^{\oplus I}$. We claim that the
map $\sigma_{L}$ factors as
$$\xymatrix{
R_{\mathscr{F}}\otimes_{R}L\ar[r]^{\simeq\:\:\:\:\:\:}
&(R_{\mathscr{F}}\otimes_{R}R)^{\oplus I}
\ar[r]^{\:\:\:\:\:\:\:\:\:(\sigma_{R})_{I}}&
(R_{\mathscr{F}})^{\oplus I}\ar[r]^{\tau^{(I)}}&(R^{\oplus I})_{\mathscr{F}}\ar[r]^{\:\:\:\:\:\:\theta_{\mathscr{F}}}&L_{\mathscr{F}}}$$
where $\tau^{(I)}\big((z_{i})_{i\in I}\big)=\sum\limits_{i\in I}(\tau_{i})_{\mathscr{F}}(z_{i})$ and for each $i$, $\tau_{i}:R\rightarrow R^{\oplus I}=\bigoplus\limits_{j\in I}R$ is the canonical map which is defined as $r\rightsquigarrow(r\delta_{i,j})_{j\in I}$. In fact, for each pure tensor $z\otimes x$ of $R_{\mathscr{F}}\otimes_{R}L$ we have $z\otimes x\rightsquigarrow(z\otimes r_{i})_{i\in I}\rightsquigarrow(r_{i}.z)_{i\in I}\rightsquigarrow\xi\rightsquigarrow\theta_{\mathscr{F}}(\xi)$ where $\xi=\sum\limits_{i\in I}(\tau_{i})_{\mathscr{F}}(r_{i}.z)$. But $\theta_{\mathscr{F}}(\xi)=z\theta_{\mathscr{F}}\big(\sum\limits_{i\in I}(\tau_{i})_{\mathscr{F}}\circ j_{R}(r_{i})\big)=z\sum\limits_{i\in I}\theta_{\mathscr{F}}\circ j_{(R^{\oplus I})}\circ\tau_{i}(r_{i})=z\sum\limits_{i\in I}j_{L}\circ\theta\circ\tau_{i}(r_{i})
=zj_{L}\Big(\theta\big(\sum\limits_{i\in I}\tau_{i}(r_{i})\big)\Big)=zj_{L}(x)=\sigma_{L}(z\otimes x)$. This establishes the claim. The map $\tau^{(I)}$ is
bijective since the G-localization functor preserves direct sums. Therefore $\sigma_{L}$ is bijective. Now, by applying the tensor functor $R_{\mathscr{F}}\otimes_{R}-$ and the G-localization functor which is, by the hypotheses, right exact to the exact sequence $\xymatrix{L'\ar[r]&L\ar[r]&M\ar[r]&0}$ where $L$ and $L'$ are free $R-$modules, and also by using a special case of the five lemma we conclude that $\sigma_{M}$ is bijective: \begin{displaymath}
\xymatrix{
R_{\mathscr{F}}\otimes_{R}L'\ar[r]^{}
\ar[d]_{\sigma_{L'}}&R_{\mathscr{F}}\otimes_{R}L
\ar[r]\ar[d]_{\sigma_{L}}&
R_{\mathscr{F}}\otimes_{R}M\ar[r]\ar[d]_{\sigma_{M}}&0\\
L'_{\mathscr{F}}\ar[r]&L_{\mathscr{F}}\ar[r]&
M_{\mathscr{F}}\ar[r]&0.}\\
\end{displaymath}
$\textbf{(v)}\Rightarrow\textbf{(iii)}:$ The G-localization functor, by Corollary \ref{G-coro 1}, is additive. In every abelian category, each split and exact sequence is left split and exact by an additive functor. This, in particular, implies that the G-localization functor preserves finite direct sums. Specially,
$\tau^{(n)}:(R_{\mathscr{F}})^{n}\rightarrow(R^{n})_{\mathscr{F}}$ is an isomorphism as $R_{\mathscr{F}}-$modules for all natural numbers $n$. Therefore all of the factors in the decomposition of $\sigma_{R^{n}}$ are bijective. Let $I$ be a f.g. ideal of $R$. Now, by applying the tensor functor $R_{\mathscr{F}}\otimes_{R}-$ and the G-localization functor which is, by the hypotheses, right exact to the exact sequence $\xymatrix{R^{n}\ar[r]&R\ar[r]&R/I\ar[r]&0}$
and also by using a special case of the five lemma we obtain that $\sigma_{R/I}$ is bijective: \begin{displaymath}
\xymatrix{
R_{\mathscr{F}}\otimes_{R}R^{n}\ar[r]^{}
\ar[d]_{\sigma_{R^{n}}}&R_{\mathscr{F}}\otimes_{R}R
\ar[r]\ar[d]_{\sigma_{R}}&
R_{\mathscr{F}}\otimes_{R}R/I\ar[r]\ar[d]_{\sigma_{R/I}}&0\\
(R^{n})_{\mathscr{F}}\ar[r]&R_{\mathscr{F}}\ar[r]&
(R/I)_{\mathscr{F}}\ar[r]&0.}
\end{displaymath}
If moreover $I\in\mathscr{F}$ then $0=(R/I)_{\mathscr{F}}\simeq R_{\mathscr{F}}\otimes_{R}R/I\simeq R_{\mathscr{F}}/IR_{\mathscr{F}}$. Therefore $IR_{\mathscr{F}}=R_{\mathscr{F}}$. In fact the latter holds for each element of the system $\mathscr{F}$ since it is of finite type.\\
$\textbf{(iii)}\Rightarrow\textbf{(vi)}:$ See Lemma \ref{prop 431}.\\
$\textbf{(vi)}\Rightarrow\textbf{(vii)}:$ Let $N$ be a $R_{\mathscr{F}}-$module. By the hypothesis, there exists a $R-$module $M$ such that $N\simeq M_{\mathscr{F}}$. Therefore $\mathscr{F}(N)\simeq\mathscr{F}(M_{\mathscr{F}})=0$.\\
$\textbf{(vii)}\Rightarrow\textbf{(iii)}:$ If $I\in\mathscr{F}$ then $R_{\mathscr{F}}\otimes_{R}R/I=0$ since  $I\subseteq\Ann_{R}\big(z\otimes(r+I)\big)$ thus $z\otimes(r+I)\in\mathscr{F}(R_{\mathscr{F}}\otimes_{R}R/I)=0$. Therefore $IR_{\mathscr{F}}=R_{\mathscr{F}}$ for all $I\in\mathscr{F}$. $\Box$ \\

\begin{proposition}\label{G-prop 2} Let $\phi:R\rightarrow S$ be a ring map, let $\mathscr{F}$ be an of finite type system on the ring $R$ and let $\mathscr{G}$ be the set of ideals $J$ of $S$ such that $S/J$ as $R-$module is $\mathscr{F}-$negligible. Then the following conditions hold.\\
$\textbf{(i)}$ An ideal $J$ of  $S$ belongs to $\mathscr{G}$ if and only if  $IS\subseteq J$ for some $I\in\mathscr{F}$.\\
$\textbf{(ii)}$ The family $\mathscr{G}$ is an idempotent topologizing system on the ring $S$ which is also of finite type. \\
$\textbf{(iii)}$ For each $S-$module $M$, $\mathscr{F}(M)=\mathscr{G}(M)$ and there is a canonical isomorphism of $R-$modules
$\eta:M_{\mathscr{G}}\rightarrow M_{\mathscr{F}}$ such that $j_{M}=\eta\circ j'_{M}$ where $j'_{M}:M\rightarrow M_{\mathscr{G}}$ is the canonical map.\\
$\textbf{(iv)}$ The map $\phi_{\mathscr{F}}:R_{\mathscr{F}}\rightarrow S_{\mathscr{F}}$ is a ring homomorphism when we identify $S_{\mathscr{F}}$ with $S_{\mathscr{G}}$ via $\eta$.\\
\end{proposition}

{\bf Proof.} $\textbf{(i)}:$ Suppose $J$ is an ideal of
$S$ and $IS\subseteq J$ for some $I\in\mathscr{F}$. Let $s\in S$ we have then
$I\subseteq\Ann_{R}(s+J)$. Therefore $S/J$ is $\mathscr{F}-$negligible. Conversely, suppose $J\in\mathscr{G}$. Then $I=\Ann_{R}(1_{S}+J)\in\mathscr{F}$ and clearly $IS\subseteq J$.\\
$\textbf{(ii)}:$ Clearly $\mathscr{G}$ is a topologizing system on the ring $S$, it is also of finite type once we have verified that it is idempotent. Let $J\in\mathscr{G}.\mathscr{G}$ then there exists an ideal $J'\in\mathscr{G}$ containing $J$ such that $J'/J$ is $\mathscr{G}-$negligible. By the hypothesis, there is a finitely generated ideal $I=\langle a_{1},...,a_{n}\rangle$ of $R$ such that $I\in\mathscr{F}$ and $IS\subseteq J'$.
For each $j$, there is also an ideal $I_{j}\in\mathscr{F}$ such that $(I_{j}a_{j})S\subseteq J$. Let $I'=I_{1}I_{2}...I_{n}$ then clearly $II'\in\mathscr{F}$ and $II'S\subseteq J$. This means that $J\in\mathscr{G}$.\\
$\textbf{(iii)}:$ If $m\in\mathscr{F}(M)$ then $I=\Ann_{R}(m)\in\mathscr{F}$ and $IS\subseteq\Ann_{S}(m)$ therefore $m\in\mathscr{G}(M)$. Conversely,
if $m\in\mathscr{G}(M)$ then there exists an ideal $I\in\mathscr{F}$ such that $IS\subseteq\Ann_{S}(m)$. Thus $I\subseteq\Ann_{R}(m)$ and so $m\in\mathscr{F}(M)$. Therefore $\mathscr{F}(M)=\mathscr{G}(M)$. \\
Let $f:J\rightarrow M/\mathscr{G}(M)$ be a $S-$linear map where $J\in\mathscr{G}$; there exists some $I\in\mathscr{F}$ such that $IS\subseteq J$. Consider the map $\eta:M_{\mathscr{G}}\rightarrow M_{\mathscr{F}}$ given by $[f]\rightsquigarrow [f\circ\phi_{|I}]$. Note that this definition is independent of choosing such $I$. We shall prove that $\eta$ is bijective. Clearly it is injective. Let $f:I\rightarrow M/\mathscr{F}(M)$ be an $R-$linear map for some $I\in\mathscr{F}$. Then consider the map $f^{\ast}:IS\rightarrow M/\mathscr{G}(M)$ which maps each $\sum\limits_{j=1}^{n}s_{j}\phi(a_{j})$ into $\sum\limits_{j=1}^{n}s_{j}f(a_{j})$. The map $f^{\ast}$ is well-defined. Because, suppose $\sum\limits_{j=1}^{n}s_{j}\phi(a_{j})=\sum\limits_{i=1}^{m}s'_{i}\phi(b_{i})$ where $a_{j},b_{i}\in I$ and $s_{j},s'_{i}\in S$. Let $f(a_{j})=m_{j}+\mathscr{F}(M)$ and $f(b_{i})=m'_{i}+\mathscr{F}(M)$ for all $i$ and $j$. Set $m=\sum\limits_{j=1}^{n}s_{j}m_{j}$ and $m'=\sum\limits_{i=1}^{m}s'_{i}m'_{i}$.
We have $IS\subseteq\Ann_{S}\big(m-m'+\mathscr{G}(M)\big)$. Because for each $c\in I$, $\phi(c)\big(m-m'+\mathscr{G}(M)\big)=
\sum\limits_{j=1}^{n}s_{j}\phi(c)f(a_{j})-
\sum\limits_{i=1}^{m}s'_{i}\phi(c)f(b_{i})=
\sum\limits_{j=1}^{n}s_{j}f(a_{j}c)-\sum\limits_{i=1}^{m}s'_{i}f(b_{i}c)=
\big(\sum\limits_{j=1}^{n}s_{j}\phi(a_{j})-
\sum\limits_{i=1}^{m}s'_{i}\phi(b_{i})\big)f(c)=0$. Therefore $m-m'+\mathscr{G}(M)\in\mathscr{G}\big(M/\mathscr{G}(M)\big)=0$. Thus $f^{\ast}$ is well-defined. Clearly it is $S-$linear and $\eta([f^{\ast}])=[f]$ and so $\eta$ is surjective.\\
$\textbf{(iv)}$: One can easily verify that $\eta^{-1}\circ\phi_{\mathscr{F}}:R_{\mathscr{F}}\rightarrow S_{\mathscr{G}}$ is actually a ring homomorphism. $\Box$ \\

\begin{lemma}\label{G-lemma 12} Let $\mathscr{F}$ be an idempotent topologizing system on the ring $R$ such that the zero ideal of $R$ does not belong to $\mathscr{F}$. If $R$ is an integral domain, then so is $R_{\mathscr{F}}$.\\
\end{lemma}

{\bf Proof.} Clearly $R_{\mathscr{F}}$ is a non-trivial ring if and only if $0\notin\mathscr{F}$. We also have $\mathscr{F}(R)=0$. Let $z,z'\in R_{\mathscr{F}}$ with representations $f$ and $g$, i.e., $z=[f]$ and $z'=[g]$ such that $z.z'=0$. Suppose $z,z'\neq0$. We may assume that the $R-$linear maps $f$ and $g$ have the same domain $I\in\mathscr{F}$, i.e.,
$f,g:I\rightarrow R$. There exists an ideal $J\in\mathscr{F}$ such that $J\subseteq f^{-1}(I)$ and $g\big(f(a)\big)=0$ for all $a\in J$. There are also elements $b,c\in J$ such that $f(b), g(c)\neq0$ since $z,z'\neq0$. But we have $0=g\big(f(bc)\big)=g\big(cf(b)\big)=f(b)g(c)$ which is a contradiction. Therefore we have either $z=0$ or that $z'=0$. $\Box$ \\

\begin{remark} If $M'$ is a $R-$submodule of $M$ then $i_{\mathscr{F}}:M'_{\mathscr{F}}\rightarrow M_{\mathscr{F}}$ is injective where $i:M'\rightarrow M$ is the canonical injection. By abuse of the notation, the image of $i_{\mathscr{F}}$ is also denoted by $M'_{\mathscr{F}}$.  \\
\end{remark}

\begin{remark}\label{remark 762} For each ideal $I$ of $R$ we have $IR_{\mathscr{F}}\subseteq I_{\mathscr{F}}$.
Moreover, $I_{\mathscr{F}}=R_{\mathscr{F}}$ for all $I\in\mathscr{F}$ because $R/I$ is $\mathscr{F}-$negligible then apply Lemma \ref{G-lemma 6}. As a second proof, we observe that $\pi|_{I}=\overline{i}\circ f$ where $\pi:R\rightarrow R/\mathscr{F}(R)$ is the canonical map, $f:I\rightarrow I/\mathscr{F}(I)$ which maps each $a\in I$ into $a+\mathscr{F}(I)$ and $i:I\rightarrow R$ is the canonical injection. Thus $1_{R_{\mathscr{F}}}=[\pi]=[\pi|_{I}]\in I_{\mathscr{F}}$. \\
\end{remark}

\begin{lemma}\label{G-lemma 13} Let $M$ be a $R-$module and let $N$ be a $R-$submodule of $M_{\mathscr{F}}$. Then $\big(j^{-1}_{M}(N)\big)_{\mathscr{F}}=
j^{-1}_{M_{\mathscr{F}}}(N_{\mathscr{F}})$. If moreover $\mathscr{F}\big(M_{\mathscr{F}}/N\big)=0$, then $N=j^{-1}_{M_{\mathscr{F}}}(N_{\mathscr{F}})$. \\
\end{lemma}

{\bf Proof.} Let $N'=j^{-1}_{M}(N)$. First we show that $N'_{\mathscr{F}}\subseteq j^{-1}_{M_{\mathscr{F}}}(N_{\mathscr{F}})$. Let $x=[\overline{i}\circ f]\in N'_{\mathscr{F}}$ where $f:I\rightarrow N'/\mathscr{F}(N')$ is an $R-$linear map and $I\in\mathscr{F}$. Consider the map $g:I\rightarrow N$ which maps each $a\in I$ into $j_{M}(z)$ where $f(a)=z+\mathscr{F}(N')$. Note that the map $g$ is well-defined because if there exists another $z'\in N'$ such that $f(a)=z'+\mathscr{F}(N')$ then $j_{M}(z-z')\in\mathscr{F}(N)\subseteq\mathscr{F}(M_{\mathscr{F}})=0$. Clearly it is also $R-$linear. Therefore $[g]\in N_{\mathscr{F}}$. We have $(\delta_{x})|_{I}=i\circ g$ where $i:N\rightarrow M_{\mathscr{F}}$ is the canonical injection. This implies that $j_{M_{\mathscr{F}}}(x)\in N_{\mathscr{F}}$. To prove the converse inclusion we act as follows. Let $x\in j^{-1}_{M_{\mathscr{F}}}(N_{\mathscr{F}})$. Thus there exists an element $y\in N_{\mathscr{F}}$ such that $j_{M_{\mathscr{F}}}(x)=
i_{\mathscr{F}}(y)$. We claim that $\Coker(\psi)=N/j_{M}(N')$ is $\mathscr{F}-$negligible where $\psi=(j_{M})|_{N'}:N'\rightarrow N$. For each $x=[f]\in N$, we show that $I\subseteq\Ann_{R}\big(x+j_{M}(N')\big)$
where $f:I\rightarrow M/\mathscr{F}(M)$ is an $R-$linear map and $I\in\mathscr{F}$. For each $a\in I$ there is some $m\in M$ such that $f(a)=\overline{m}=m+\mathscr{F}(M)$. Clearly $(\delta_{\overline{m}})|_{I}=a.f$ thus $j_{M}(m)=a.x\in N$ and so $m\in N'$.
This establishes the claim. We also have $\Ker(\psi)=\mathscr{F}(N')$. Therefore, by Theorem \ref{G-lemma 10}, $\psi_{\mathscr{F}}$ is bijective. Thus there is an element $z\in N'_{\mathscr{F}}$ such that $\psi_{\mathscr{F}}(z)=y$. Therefore $x=i'_{\mathscr{F}}(z)\in\Ima i'_{\mathscr{F}}=N'_{\mathscr{F}}$ because $j_{M_{\mathscr{F}}}\circ i'_{\mathscr{F}}=i_{\mathscr{F}}\circ\psi_{\mathscr{F}}$ where $i':N'\rightarrow M$ is the canonical injection. Now we prove the last part of the assertion. For each $x\in N$, $\delta_{x}=i\circ g$ where the map $g:R\rightarrow N$ is defined as $r\rightsquigarrow r.x$. Thus $N\subseteq j^{-1}_{M_{\mathscr{F}}}(N_{\mathscr{F}})$. Conversely, let
$x\in j^{-1}_{M_{\mathscr{F}}}(N_{\mathscr{F}})$ then there exists an ideal $I\in\mathscr{F}$ such that $I\subseteq\Ann_{R}(x+N)$. Therefore $x+N\in\mathscr{F}\big(M_{\mathscr{F}}/N\big)=0$. $\Box$ \\

Now we prove the final main result of this article:\\

\begin{theorem}\label{theorem 33} Let $\mathscr{F}$ be an of finite type system on the ring $R$ and let $\mathscr{F}'$ be the set of ideals $J$ of $R_{\mathscr{F}}$ such that $R_{\mathscr{F}}/J$ as $R-$module is $\mathscr{F}-$negligible. Then the map $\mathfrak{p}\rightsquigarrow\mathfrak{p}_{\mathscr{F}}$ is a bijection between the set of prime ideals of $R$ which do not belong to $\mathscr{F}$ and the set of prime ideals of $R_{\mathscr{F}}$ which do not belong to $\mathscr{F}'$.\\
\end{theorem}

{\bf Proof.} First of all we show that if $\mathfrak{p}$ is a prime ideal of $R$ such that $\mathfrak{p}\notin\mathscr{F}$ then
$\mathfrak{p}_{\mathscr{F}}$ is a prime ideal of $R_{\mathscr{F}}$ and that $\mathfrak{p}_{\mathscr{F}}\notin\mathscr{F}'$. Suppose $1_{R_{\mathscr{F}}}=[\pi']\in\mathfrak{p}_{\mathscr{F}}$ where $\pi':R\rightarrow R/\mathscr{F}(R)$ is the canonical map. Thus there exists an ideal $I\in\mathscr{F}$ such that for each $a\in I$ there is some $x_{a}\in\mathfrak{p}$ in which $\Ann_{R}(a-x_{a})\in\mathscr{F}$. But $I$ is not contained in $\mathfrak{p}$ since $\mathfrak{p}\notin\mathscr{F}$. Take $b\in I\setminus\mathfrak{p}$ then clearly $\Ann_{R}(b-x_{b})\subseteq\mathfrak{p}$. This implies that $\mathfrak{p}\in\mathscr{F}$ which is a contradiction therefore $\mathfrak{p}_{\mathscr{F}}$ is a proper ideal of $R_{\mathscr{F}}$.\\
We shall apply Proposition \ref{G-prop 2} to the canonical ring map $\pi:R\rightarrow R/\mathfrak{p}$. If $0\in\mathscr{G}$ then $R/\mathfrak{p}$ will be $\mathscr{F}-$negligible. By applying Lemma \ref{G-lemma 6} to the exact sequence $\xymatrix{0\ar[r]&\mathfrak{p}\ar[r]&R\ar[r]^{\pi\:\:\:\:}&R/\mathfrak{p}}$
we get that $\mathfrak{p}_{\mathscr{F}}=R_{\mathscr{F}}$ which is a contradiction since $\mathfrak{p}_{\mathscr{F}}$ is a proper ideal of $R_{\mathscr{F}}$. Thus $0\notin\mathscr{G}$ and so by Lemma \ref{G-lemma 12}, $(R/\mathfrak{p})_{\mathscr{G}}\simeq(R/\mathfrak{p})_{\mathscr{F}}$ is an integral domain. Suppose $xy\in\mathfrak{p}_{\mathscr{F}}$ for some elements $x,y\in R_{\mathscr{F}}$. Thus $0=\pi_{\mathscr{F}}(xy)=\pi_{\mathscr{F}}(x)\pi_{\mathscr{F}}(y)$. Therefore we have either $\pi_{\mathscr{F}}(x)=0$ or that $\pi_{\mathscr{F}}(y)=0$. Hence $\mathfrak{p}_{\mathscr{F}}$ is a prime ideal.\\
Suppose $\mathfrak{p}_{\mathscr{F}}\in\mathscr{F}'$ then there exists an ideal $I\in\mathscr{F}$ such that $IR_{\mathscr{F}}\subseteq\mathfrak{p}_{\mathscr{F}}$. This implies that $I\subseteq\mathfrak{p}$. Because for each $a\in I$ there exists an ideal $J\in\mathscr{F}$ such that for each $b\in J$ there is some $y_{b}\in\mathfrak{p}$ for which $\Ann_{R}(ab-y_{b})\in\mathscr{F}$. Take some $b\in J\setminus\mathfrak{p}$ then $\Ann_{R}(ab-y_{b})$ is not contained in $\mathfrak{p}$. It follows that $a\in\mathfrak{p}$. But $I\subseteq\mathfrak{p}$ is a contradiction since $\mathfrak{p}\notin\mathscr{F}$, therefore $\mathfrak{p}_{\mathscr{F}}\notin\mathscr{F}'$. \\ The map $\mathfrak{p}\rightsquigarrow\mathfrak{p}_{\mathscr{F}}$ is injective between the foregoing sets. Because suppose $\mathfrak{p}_{\mathscr{F}}=\mathfrak{p}'_{\mathscr{F}}$ where $\mathfrak{p}$ and $\mathfrak{p}'$ are prime ideals of $R$ with $\mathfrak{p},\mathfrak{p}'\notin\mathscr{F}$. For each $a\in\mathfrak{p}$, $j_{\mathfrak{p}}(a)\in\mathfrak{p}'_{\mathscr{F}}$ thus there exists an ideal $J\in\mathscr{F}$ such that for each $b\in J$ there is some $x_{b}\in\mathfrak{p}'$ in which $\Ann_{R}(ab-x_{b})\in\mathscr{F}$. Now take some $b\in J\setminus\mathfrak{p}'$ also take some $c\in\Ann_{R}(ab-x_{b})\setminus\mathfrak{p}'$, then we have $c(ab-x_{b})=0\in\mathfrak{p}'$ and so $a\in\mathfrak{p}'$; symmetrically we have $\mathfrak{p}'\subseteq\mathfrak{p}$, thus $\mathfrak{p}=\mathfrak{p}'$. \\
Finally, we show that the map $\mathfrak{p}\rightsquigarrow\mathfrak{p}_{\mathscr{F}}$  is surjective. Let $\mathfrak{q}$ be a prime ideal of $R_{\mathscr{F}}$ which does not belong to $\mathscr{F}'$. We have $\mathfrak{p}=j_{R}^{-1}(\mathfrak{q})\notin\mathscr{F}$ since $\mathfrak{p}R_{\mathscr{F}}\subseteq\mathfrak{q}$. We also have  $\mathscr{F}(R_{\mathscr{F}}/\mathfrak{q})=0$ because pick $x+\mathfrak{q}\in\mathscr{F}(R_{\mathscr{F}}/\mathfrak{q})$, if $x\notin\mathfrak{q}$ then $\Ann_{R}(x+\mathfrak{q})R_{\mathscr{F}}\subseteq\mathfrak{q}$ it follows that $\mathfrak{q}\in\mathscr{F}'$ which is a contradiction.
Therefore, by Lemma \ref{G-lemma 13}, $\mathfrak{p}_{\mathscr{F}}=\mathfrak{q}$. $\Box$ \\

The usual localization theory is just a special case of the G-localizations:\\

\begin{proposition}\label{G-th 1} Let $S$ be a multiplicative subset in the ring $R$ and let $\mathscr{F}$ be the set of ideals of $R$ which meeting $S$. Then for each $R-$module $M$, there is a unique $R-$linear map $\phi: S^{-1}M\rightarrow M_{\mathscr{F}}$ such that $j_{M}=\phi\circ\pi$ where $\pi:M\rightarrow S^{-1}M$ is the canonical map. Moreover $\phi$ is bijective.
In particular, the $R-$algebras $S^{-1}R$ and $R_{\mathscr{F}}$ are canonically isomorphic.\\
\end{proposition}

{\bf Proof.} Clearly $\Ker(\pi)$ and $\Coker(\pi)$ are $\mathscr{F}-$negligible. Therefore, by Proposition \ref{G-lemma 2}, the first part of the assertion is realized. Now for each $I\in\mathscr{F}$, consider the map $\lambda_{I}:\Hom_{R}\big(I,M/\mathscr{F}(M)\big)\rightarrow S^{-1}M$ defined by $f\rightsquigarrow m/s$ where $f(s)=m+\mathscr{F}(M)$ for some $s\in I\cap S$. The map $\lambda_{I}$ is well-defined. Because suppose there is another element $s'\in I\cap S$ such that $f(s')=m'+\mathscr{F}(M)$. We have $s'f(s)=f(ss')=sf(s')$. It implies that $S\cap\Ann_{R}(sm'-s'm)\neq\emptyset$. The map $\lambda_{I}$ is also $R-$linear and clearly $\lambda_{I}=\lambda_{J}\circ u_{I,J}$ for every $I,J\in\mathscr{F}$ with $I\leq J$. Therefore, by the universal property of the colimits, there is a (unique) $R-$linear map $\lambda:M_{\mathscr{F}}\rightarrow S^{-1}M$ such that the map $\lambda_{I}$ factors as
$\xymatrix{\Hom_{R}\big(I,M/\mathscr{F}(M)\big)\ar[r]
& M_{\mathscr{F}}\ar[r]^{\lambda\:\:\:\:\:}&S^{-1}M}$ for all $I\in\mathscr{F}$. The map $\lambda$ is injective because suppose $\lambda([f])=0$ for some element $[f]\in M_{\mathscr{F}}$ where $f:I\rightarrow M/\mathscr{F}(M)$ is a $R-$linear map and $I\in\mathscr{F}$. We know that $I\cap S\neq\emptyset$. Take some $s\in I\cap S$ then clearly $Rs\in\mathscr{F}$ and $f|_{Rs}=0$. The map $\lambda$ is also surjective. Because for each element $m/s\in S^{-1}M$ then consider the $R-$linear map $g:Rs\rightarrow M/\mathscr{F}(M)$ defined by $rs\rightsquigarrow rm+\mathscr{F}(M)$. Note that the map $g$ is well-defined since if $rs=r's$ then $r-r'\in\mathscr{F}(R)$ and so $(r-r')m\in\mathscr{F}(R)M\subseteq\mathscr{F}(M)$. It is also $R-$linear and clearly $\lambda([g])=m/s$. We have $\lambda^{-1}=\phi$. $\Box$ \\

\end{document}